\documentclass{svproc}
\usepackage{url}

\usepackage{hyperref}
\usepackage{type1cm}        
\usepackage{makeidx}         
\usepackage{graphicx}        

\usepackage{multicol}        
\usepackage[bottom]{footmisc}

\usepackage{newtxtext}       %
\usepackage[varvw]{newtxmath}       

\usepackage{bm}
\usepackage{siunitx}

\newcounter{algorithm}

\newcommand{\bbC}{{\mathbb{C}}}

\newcommand{\bbE}{{\mathbb{E}}}

\newcommand{\bbI}{{\mathbb{I}}}

\newcommand{\bbN}{{\mathbb{N}}}

\newcommand{\bbP}{{\mathbb{P}}}

\newcommand{\bbR}{{\mathbb{R}}}

\DeclareSymbolFont{bbold}{U}{bbold}{m}{n}
\DeclareSymbolFontAlphabet{\mathbbold}{bbold}

\providecommand{\argmin}{\operatorname*{argmin}}

\usepackage{sidecap}

\usepackage{color}

\def\D{\,\mathrm{d}}

\newcommand*{\supp}{\mathrm{supp}}

\newcommand*{\ord}[1]{\cO ( #1  ) }

\newcommand*{\nm}[1]{{\|#1\|}}

\DeclareMathOperator*{\esssup}{ess\,sup}

\newcommand{\tnm}[1]{{\left\vert\kern-0.25ex\left\vert\kern-0.25ex\left\vert #1 
    \right\vert\kern-0.25ex\right\vert\kern-0.25ex\right\vert}}

\def\ef#1{\eqref{#1}}

\DefineNamedColor{named}{Purple}{cmyk}{0.45,0.86,0,0}

\DefineNamedColor{named}{JungleGreen} {cmyk}{0.99,0,0.52,0}

\DefineNamedColor{named}{orange} {rgb}{1,0.55,0}

\newcommand*{\cD}{\mathcal{D}}

\newcommand*{\cK}{\mathcal{K}}

\newcommand*{\cR}{\mathcal{R}}
\newcommand*{\cO}{\mathcal{O}}
\newcommand*{\cP}{\mathcal{P}}

\newcommand{\be}[1]{\begin{equation} #1 \end{equation}}
\newcommand{\bes}[1]{\begin{equation*} #1 \end{equation*}}

\newcommand{\eas}[1]{\begin{align*} #1 \end{align*}}

\newcommand{\thm}[1]{\begin{theorem} #1 \end{theorem}}

\newcommand{\rem}[1]{\begin{remark} #1 \end{remark}}
\newcommand{\lem}[1]{\begin{lemma} #1 \end{lemma}}
\newcommand{\cor}[1]{\begin{corollary} #1 \end{corollary}}

\newcommand{\prf}[1]{\begin{proof} #1 \end{proof}}

\begin{document}

\mainmatter              
\title{Function recovery and optimal sampling in the presence of nonuniform evaluation costs}
\titlerunning{Function recovery with nonuniform costs}  
%
\author{Ben Adcock\inst{1}}
\authorrunning{Ben Adcock} 
%
\tocauthor{Ben Adcock}
\institute{Simon Fraser University, Burnaby, BC, Canada\\
\email{ben\_adcock@sfu.ca}}


\maketitle              

\begin{abstract}
We consider recovering a function $f : D \rightarrow \bbC$ in an $n$-dimensional linear subspace $\cP$ from i.i.d.\ pointwise samples via (weighted) least-squares estimators. Different from most works, we assume the cost of evaluating $f$ is potentially nonuniform, and governed by a cost function $c : D \rightarrow (0,\infty)$ which may blow up at certain points. We therefore strive to choose the sampling measure in a way that minimizes the expected total cost.  We provide a recovery guarantee which asserts accurate and stable recovery with an expected cost depending on the \textit{Christoffel function} and \textit{Remez constant} of the space $\cP$. This leads to a general recipe for finding a good sampling measure for general $c$. As an example, we consider one-dimensional polynomial spaces. Here, we provide two strategies for choosing the sampling measure, which we prove are optimal (up to constants and log factors) in the case of algebraically-growing cost functions. 
\end{abstract}

\section{Introduction}

Countless applications require one to construct an accurate estimator of an unknown function $f : D \rightarrow \bbC$  defined over a domain $D$ (typically $D \subseteq \bbR^d$) from noisy samples
\be{
\label{f_meas}
y_i = f(x_i) + e_i,\quad i = 1,\ldots,m.
}
In many such applications, data is expensive to acquire, yet there is significant freedom in choosing the sample points $x_i$. The construction of good points has spurred a vast amount of research across different fields, such as optimal design of experiments in statistics, active learning in machine learning, and optimal information in information-based complexity. 

In recent years, this question has been essentially settled in the case where one strives to construct an estimator $\hat{f}$ in a fixed, but arbitrary $n$-dimensional subspace $\cP \subseteq L^2_{\varrho}(D)$. Drawing the samples $x_1,\ldots,x_m$ i.i.d.\ according to the probability measure $\mu$ whose density  is proportional to the \textit{Christoffel function} of $\cP$ (so-called \textit{Christoffel sampling}) is a near-optimal random sampling strategy \cite{cohen2017optimal,hampton2015coherence} (see \cite{adcock2024optimala} for a review). Specifically, the estimator $\hat{f}$ obtained via a weighted least-squares fit to the data \ef{f_meas} yields a near-best recovery of $f$ from $\cP$ provided the number of samples  $m \gtrsim n \log(n)$ scales log-linearly in $n = \dim(\cP)$. Refinements of this approach can also further reduce the sampling budget by removing the $\log(n)$ term (see \cite{bartel2023constructive,dolbeault2024randomized,krieg2021function,krieg2021functionII} and references therein).

\subsection{Optimal sampling with nonuniform evaluation costs}

It is common in applications to assume a fixed cost-per-evaluation model, i.e., evaluating $f(x_i)$ involves the same amount of resources (computational, physical, or otherwise) regardless of the choice of $x_i$. However, in practice evaluation costs may be highly variable depending on the location of $x_i$.

Nonuniform evaluation costs are problematic when any sampling strategy optimized for a uniform cost model strives to place more samples at locations with higher evaluation costs. In the case of Christoffel sampling, for instance, this in effect means regions where the Christoffel function is large. While Christoffel sampling is optimal (up to constants and log factors) for reducing the total number of samples, it may be highly suboptimal in terms of the total evaluation cost.

With this in mind, in this paper we consider function recovery in the setting of nonuniform evaluation costs and, specifically, how this affects the choice of sampling measure. We consider a \textit{cost function} $c : D \rightarrow (0,\infty)$ and, rather than the problem of optimizing the $x_i$ for a fixed maximum sampling budget $m$, we consider the problem of optimizing the $x_i$ subject to a fixed total cost $C_{\mathsf{tot}} = c(x_1) + \cdots + c(x_m)$. We focus on weighted least-squares estimators in arbitrary $n$-dimensional subspaces $\cP$ with sample points $x_1,\ldots,x_m$ chosen i.i.d.\ with respect to some probability measure $\mu$. Thus, the matter at hand is how to chose $\mu$ judiciously in terms of $\cP$ \textit{and} $c$ to ensure as accurate an approximation to $f$ as possible.

\subsection{Setup}\label{ss:setup}

We now describe our setup in more detail. This is based on \cite{adcock2024optimala}. Let $(D,\cD,\varrho)$ be a probability space and $L^2_{\varrho}(D)$ be the Lebesgue space of square-integrable functions $f : D \rightarrow \bbC$ with respect to $\varrho$. 
Suppose that $\mu$ that is absolutely continuous with respect to $\varrho$, let $v = \D \mu / \D \varrho$ and define the weight function  $w : D \rightarrow [0,\infty]$ as 
\be{
\label{weight-fn}
w(x) = 1/v(x),\quad x \in D.
}
Let $f \in L^2_{\varrho}(D)$ be the function to recover, $x_1,\ldots,x_m \sim_{\mathrm{i.i.d.}} \mu$ and consider samples \ef{f_meas}. Let $\cP \subset L^2_{\varrho}(D)$ be an arbitrary $n$-dimensional subspace, termed the \textit{approximation space}. We make the mild, but convenient assumption that $1 \in \cP$, where $1$ is the constant function equal to one almost everywhere. 
We then consider 
weighted least-squares fit
\be{
\label{wls-prob}
\min_{p \in \cP} \frac1m \sum^{m}_{i=1} w(x_i) | y_i - p(x_i) |^2.
}
Despite involving $L^2$ functions, this problem is well-defined with probability one, since the sample points $x_i$ are drawn randomly from $\mu \ll \varrho$.
Let $c$ and $C_{\mathsf{tot}}$ be as above. 
Since $C_{\mathsf{tot}}$ is a random variable, we define the \textit{expected cost}
\be{
\label{C-exp-def}
C_{\mathsf{exp}} = \bbE(C_{\mathsf{tot}}) = m \cdot \int_{D} c(x) \D \mu(x)
}
(we write $C_{\mathsf{exp}} = + \infty$ if $c$ is not integrable with respect to $\mu$).
In this paper, we are interested in conditions that ensure the approximation $\hat{f}$ given by \ef{wls-prob} is \textit{quasi-best} in the sense that the $L^2_{\varrho}$-norm error $\nm{f - \hat{f}}_{L^2_{\varrho}(D)}$ is proportional (up to the noise $e_i$) to the \textit{best approximation error}
\bes{
\inf_{p \in \cP} \tnm{f - p} 
}
measured in some suitable norm $\tnm{\cdot}$. With this in hand, the main questions we address in this paper are:
\begin{enumerate}
\item[(i)] Given a measure $\mu$ and cost function $c$, how large should the expected cost $C_{\mathsf{exp}}$ be to ensure a quasi-best approximation?
\item[(ii)] How should we choose $\mu$ to minimize $C_{\mathsf{exp}}$?
\end{enumerate}
We remark in passing that the noise terms $e_i$ in \ef{f_meas} are not the focus of this paper. For concreteness, we assume the noise is deterministic and bounded, and aim to show that the error of the estimator $\hat{f}$ scales linearly in $\nm{e}_2$, where $e = (e_i)^{m}_{i=1}$. Other noise models, including statistical models, could readily be considered instead.

\subsection{Motivations and related work}

This work is motivated by surrogate model construction in UQ. Here, $f$ is typically a (quantity-of-interest) of a parametric PDE, with $x \in D$ representing the parameters. A black-box numerical simulation is used to evaluate $f(x)$, which often involves significant computational cost. It is common to assume constant cost for any $x \in D$. However, in practice, the cost may vary significantly with $x$. For instance, certain parameter values may lead to PDE solutions involving shocks, oscillations, boundary layers or other behaviours, all of which require finer meshes in any (adaptive) PDE discretization.

This problem has commonly been tackled by \textit{multi-level} and \textit{multi-fidelity} methods \cite{giles2015multilevel,peherstorfer2018survey}. In these approaches, one is able to access approximations $f_h(x) \approx f(x)$ with different accuracies (fidelities) and different evaluation costs, which are then combined to design accurate estimators given a fixed cost budget. In this work, we take a different perspective. We assume $f$ is evaluated by a single model, represented as a true black-box, where the cost of an evaluation is determined solely by the only input $x$. Our work is also inspired by \cite{oates2024probabilistic}, although it addresses a different problem. When we specialize to polynomials, our approach (which uses extrapolation estimates) is also related to \cite{demanet2018stable}.

Finally, we remark in passing that nonuniform evaluation costs (or \textit{labelling} costs, as they are more commonly known) arise in various machine learning applications \cite{settles2008active}. Various \textit{cost-sensitive} learning strategies have been developed to handle this setting.

\subsection{Contributions and outline}

We commence in \S \ref{s:background} with background, focusing on Christoffel functions and Christoffel sampling. We also introduce our main example, which involves computing polynomial estimators over the domain $D = (-1,1)$.
In \S \ref{s:recov-general} we consider recovery in the general setting (arbitrary $D$, $\varrho$ and $\cP$). As noted, the issue with nonuniform evaluation costs is that large regions of the cost function may correspond to regions one may want to sample more densely to ensure an accurate estimator. With this in mind, we establish a result, Theorem \ref{thm:accuracy-fulldomain}, which blends accurate recovery in subdomains where the cost function is small (and dense sampling is therefore permitted) with extrapolation to regions where it is large and sample points are scarce. This yields a general recipe for constructing sampling measures for nonuniform evaluation costs, given in terms of the Christoffel function over a subdomain $\Omega$, which determines the accuracy of recovery in $\Omega$, and the so-called \textit{Remez constant}, which ensures accurate extrapolation outside of $\Omega$. 

In \S \ref{s:recov-poly} we consider the main example. We introduce two approaches for constructing $\mu$ in this case. The first assumes knowledge of the cost function $c$, but is inherently \textit{hierarchical}, since the resulting measure is independent of $n = \dim(\cP)$. The second approach is nonhierarchical, but requires no knowledge of $c$. When the cost function is of the form $c(x) = (1-x^2)^{-\alpha}$ for $\alpha \geq 1$, the expected cost $C_{\mathsf{exp}}$ obtained through this approach  behaves like $n^{2 \alpha} \log(n)$. As we discuss, this scaling is effectively optimal.
We end \S \ref{s:recov-poly} with some numerical experiments and then offer some conclusions in \S \ref{s:conclusions}.

\section{Background and main example}\label{s:background}

We first require some background. For a finite-dimensional subspace $\cP \subseteq L^2_{\varrho}(D)$, we define its (reciprocal) \textit{Christoffel function} as 
\be{
\label{Kappa-def}
\cK(x) = \cK(\cP , L^2_{\varrho}(D))(x) : = \sup  \{ | p(x) |^2 : p \in \cP,\ \nm{p}^2_{L^2_{\varrho}(D)} = 1 \},\quad \forall x \in D.
}
We remark that $\cK$ has the equivalent expression
\be{
\label{Kappa-def-alt}
\cK(x) = \sum^{n}_{i=1} | \phi_i(x) |^2,
}
where $\{ \phi_i \}^{n}_{i=1}$ is any orthonormal basis of $\cP$. Given a weight function $w$ as in \ef{weight-fn}, for convenience we also define
\be{
\label{kappa_w_def}
\kappa_w = \kappa_w(\cP , L^2_{\varrho}(D)) = \esssup_{x \sim \varrho} \left \{ w(x) \cK(\cP , L^2_{\varrho}(D) )(x) \right \}.
}
Finally, given $f \in L^2_{\varrho}(D)$, we write $f_n \in \cP$ for its best $L^2_{\varrho}(D)$-norm approximation in $\cP$, i.e.,
\be{
\label{orth-proj}
f_n = \argmin_{p \in \cP} \nm{f - p}_{L^2_{\varrho}(D)}.
}

\subsection{Christoffel functions and sample complexity}

We now present a standard result that shows how the Christoffel function relates to the sample complexity of the weighted least-squares fit \ef{wls-prob} with random sampling. Here and elsewhere, we write $\mathrm{supp}(\cdot)$ for the \textit{support} of a measure.

\thm{
\label{t:err-bd-prob}
Suppose that $v = \D \mu / \D \varrho$ is strictly positive almost everywhere on $\mathrm{\supp}(\varrho)$. Let $f \in L^2_{\varrho}(D)$, $0 < \epsilon < 1$, $\kappa = \kappa_1(\cP , L^2_{\varrho}(D)) = \esssup_{x \sim \varrho} \cK(\cP , L^2_{\varrho}(D) )(x) $ and
\be{
\label{samp-comp-standard}
m \geq 8 \cdot \kappa_w(\cP , L^2_{\varrho}(D) ) \cdot \log(3 n / \epsilon),
}
and consider $x_1,\ldots,x_m \sim_{\mathrm{i.i.d.}} \mu$.
Then the following holds with probability at least $1-\epsilon$. For any $e \in \bbC^m$, the weighted least-squares estimator $\hat{f}$ is unique and satisfies
\eas{
\nm{f - \hat{f}}_{L^2_{\varrho}(D)}  &\leq \frac{9}{\sqrt{\epsilon}} \nm{f - f_n}_{L^2_{\varrho}(D)} + \frac{8}{\sqrt{\epsilon}} \nm{e}_{2}
\\
\nm{f - \hat{f}}_{L^{\infty}_{{\varrho}}(D)} &\leq \nm{f - f_n}_{L^{\infty}_{{\varrho}}(D)} + 8\sqrt{\frac{\kappa}{\epsilon}} \nm{f-f_n}_{L^2_{\varrho}(D)} + 8\sqrt{\frac{\kappa}{\epsilon}} \nm{e}_2,
}
where $f_n$ is as in \ef{orth-proj}.
}

This theorem is a simple consequence of a result we prove later, Lemma \ref{lem:accuracy-subdomain}. Note that the constant $8$ in  \ef{samp-comp-standard} can be changed, at the expense of different constants in the error bounds. There are many variations on this result in the literature, including bounds in expectation instead of probability and `in probability' bounds that are uniform in $f$. See \cite{adcock2024optimala} for a review of these different bounds. The above formulation is particularly convenient for our purposes when we deal with nonuniform evaluation costs. However, we anticipate such variations could also be considered. 

\subsection{Christoffel sampling}\label{ss:CS}

If the evaluation cost is uniform, this theorem immediately implies an optimal random sampling strategy. One chooses $\mu$ (or, equivalently, $w$) to optimize \ef{samp-comp-standard} amongst all probability measures, resulting in the choice
\be{
\label{mu-opt}
\D \mu(x) = n^{-1} \cK(x) \D \varrho(x).
}
Note that the $n^{-1}$ terms comes from the fact that $\int_{D} \cK(x) \D \varrho(x) = n$, which is itself an immediate consequence of \ef{Kappa-def-alt}. With this choice, \ef{samp-comp-standard} reduces to $m \geq 8 n \cdot \log(3n/\epsilon)$.
Hence, we obtain optimal sample complexity up to the constant and log factor.

\subsection{Main example and intuition}\label{ss:main-examp}

In our main example, we consider the case where $D = (-1,1)$, $\varrho$ is the uniform probability measure on $D$, i.e., $\D \varrho(x) = \frac12 \D x$, and $\cP = \bbP_{n-1}$ is the space of algebraic polynomials of degree less than $n$. Polynomial estimators are ubiquitous in applications, and served as some of the primary motivations in the development of Christoffel sampling techniques \cite{adcock2024optimala,cohen2013stability,cohen2017optimal,hampton2015coherence}. In this work we consider only the one-dimensional setting for convenience. We anticipate that our main results will generalize to higher dimensions.

This space has an orthonormal basis $\{ \phi_1,\ldots,\phi_n \}$, where $\phi_{i} \in \bbP_{i-1}$ is the $(i-1)$th orthonormal Legendre polynomial (see, e.g., \cite[\S 2.2.2]{adcock2022sparse}). These polynomials satisfy
\be{
\label{leg-poly-max}
\max_{-1 \leq x \leq 1} | \phi_i(x) | = \phi_i(1) = \sqrt{2i-1},\quad \forall i \in \bbN.
}
They also satisfy the envelope bound (see, e.g., \cite[Rem.\ 5.21]{adcock2022sparse})
\be{
\label{envelope}
| \phi_i(x) | \leq \frac{\sqrt{2}}{\sqrt{\pi/2}(1-x^2)^{1/4}},\quad \forall x \in (-1,1),\ i \in \bbN.
}
As a result, \ef{Kappa-def-alt} implies that the Christoffel function of $\cP$ \ef{Kappa-def} satisfies
\be{
\label{Christoffel-UB}
\cK(x) \leq \min \left \{ n^2 , \frac{4 n}{\pi (1-x^2)^{1/2}} \right \},\qquad \sup_{x \in [-1,1]} \cK(x) = \cK(1) = n^2.
}
Moreover, one can show that the sequence of Christoffel sampling measures \ef{mu-opt} converges weakly as $n \rightarrow \infty$ to the Chebyshev (arcsine) measure
\be{
\label{arcsine-meas}
\frac{1}{\pi \sqrt{1-x^2}} \D x,
}
which is the equilibrium measure for $D = (-1,1)$ (see, e.g., \cite{narayan2017christoffel}).

This discussion illustrates the potential problems of Christoffel sampling in the presence of nonuniform evaluation costs. Christoffel sampling will place more samples near the interval endpoints $x = \pm 1$. Hence, if the cost function also blows up as $|x| \rightarrow 1^{-}$, then Christoffel sampling, despite being near-optimal in terms of the number of samples, will lead to a large evaluation cost. Indeed, the expected cost \ef{C-exp-def} may be infinite if $c$ is not integrable with respect to \ef{mu-opt} or \ef{arcsine-meas}. On the other hand, in the more benign scenario where the cost function is bounded near $x = \pm 1$, then Christoffel sampling may also be near-optimal in terms of the total evaluation cost.

\section{Recovery with arbitrary least-squares estimators}\label{s:recov-general}

With this discussion in mind, we now present some general analysis for least-squares estimators in the presence of nonuniform evaluation costs. Our approach is based on the notion of \textit{extrapolation}. Our first result, Lemma \ref{lem:accuracy-subdomain}, shows that the least-squares estimator can always provide accurate recovery over a subdomain where the samples are sufficiently dense. We then combine this in Theorem \ref{thm:accuracy-fulldomain} with extrapolation tools to provide error bounds over the full domain.

In what follows, we write $\varrho |_{\Omega}$ for the restriction of $\varrho$ to a $\varrho$-measurable set $\Omega \subseteq D$. To avoid unnecessary notation, we write $L^2_{\varrho}(\Omega)$ instead of $L^2_{\varrho |_{\Omega}}(\Omega)$ for the corresponding space of square-integrable functions. 

\subsection{Recovery in $\Omega$}

\lem{
\label{lem:accuracy-subdomain}
Let $\Omega \subseteq D$ be a $\varrho$-measurable set and suppose that $v = \D \mu / \D \varrho$ is strictly positive almost everywhere on $\Omega$.  Let $f \in L^2_{\varrho}(D)$, $0 < \epsilon < 1$, $\kappa = \kappa_1(\cP , L^2_{\varrho}(\Omega) ) = \esssup_{x \sim \varrho|_{\Omega}}  \cK(\cP , L^2_{\varrho}(\Omega) )(x) $ and
\be{
\label{m-cond-subdomain}
m \geq 8 \cdot \kappa_w(\cP , L^2_{\varrho}(\Omega) ) \cdot \log(3 n/\epsilon),
}
and consider $x_1,\ldots,x_m \sim_{\mathrm{i.i.d}} \mu$.
Then the following holds with probability at least $1-\epsilon$. For any $e \in \bbC^m$ and any minimizer $\hat{f}$ of \ef{wls-prob}, the error satisfies
\eas{
\nm{f - \hat{f}}_{L^2_{\varrho}(\Omega)}  & \leq \frac{9}{\sqrt{\epsilon}}\nm{f - f_n}_{L^2_{\varrho}(D)} + \frac{8}{\sqrt{\epsilon}} \nm{e}_{2}
\\
\nm{f - \hat{f}}_{L^{\infty}_{\varrho}(\Omega)} & \leq \nm{f - f_n}_{L^{\infty}_{\varrho}(\Omega)} + 8\sqrt{\frac{\kappa}{\epsilon}}  \nm{f-f_n}_{L^2_{\varrho}(D)} + 8\sqrt{\frac{\kappa}{\epsilon}} \nm{e}_2,
}
where $f_n$ is as in \ef{orth-proj}.
}

This lemma states that accurate recovery can occur in a subdomain with a total number of samples depending on the Christoffel function over the subdomain only. In particular, if $\Omega$ is a region where $c(x)$ is small then we can always ensure accurate recovery in this region while avoiding a high evaluation cost. In contrast with Theorem \ef{t:err-bd-prob}, the estimator $\hat{f}$ may no longer be unique. This stems from the fact that $v = \D \mu / \D \varrho$ may not be strictly positive on $\mathrm{supp}(\varrho)$. However, the same error bound holds regardless of the chosen minimizer. In practice, one might choose the minimizer with minimal $L^2_{\varrho}(D)$-norm, which can be computed via the pseudoinverse of the matrix of the algebraic least-squares problem associated to \ef{wls-prob}.

\prf{
Notice that the least-squares problem \ef{wls-prob} is almost surely well-defined, due to the assumptions on $f$ and $\mu$. A solution therefore exists, since $\cP$ is finite-dimensional. We now establish the error bounds. In the standard analysis of least-squares problems (see, e.g., \cite[Chpt.\ 5]{adcock2022sparse}), the so-called \textit{discrete stability constant} 
\bes{
\alpha_0 =  \inf \left \{ \sqrt{\frac1m \sum^{m}_{i=1} w(x_i) | p(x_i) |^2 } : p \in \cP,\ \nm{p}_{L^2_{\varrho}(D)} = 1 \right \}.
}
plays a crucial role. If $\alpha_0 > 0$ the solution is unique and it obeys an error bound over $D$ where the constants scale like $1/\alpha_0$. Since in the current lemma we consider the error over the subset $\Omega$, we now introduce the following constant instead:
\eas{
\alpha & = \inf \left \{ \sqrt{\frac1m \sum^{m}_{i=1} w(x_i) | p(x_i) |^2 \bbI_{\Omega}(x_i)}  : p \in \cP,\ \nm{p}_{L^2_{\varrho}(\Omega)} = 1 \right \} \leq \alpha_0 .
}
Write
\be{
\label{32C-today}
\nm{f - \hat{f}}_{L^2_{\varrho}(\Omega)}  \leq \nm{f - f_n}_{L^2_{\varrho}(\Omega)} + \nm{f_n - \hat{f}}_{L^2_{\varrho}(\Omega)} .
}
We now estimate the second term. Let $E_1$ be the event that $\alpha \geq 1/2$ (we will estimate the probability of this event later) and suppose that $E_1$ occurs. Then we have
\eas{
\nm{f_n - \hat{f}}_{L^2_{\varrho}(\Omega) } & \leq  2 \sqrt{\frac1m \sum^{m}_{i=1} w(x_i) | f_n(x_i) - \hat{f}(x_i) |^2 \bbI_{\Omega}(x_i)}
\\
& \leq  2 \sqrt{\frac1m \sum^{m}_{i=1} w(x_i) | y_i - \hat{f}(x_i) |^2 \bbI_{\Omega}(x_i)} 
\\
&~~~ + 2 \sqrt{\frac1m \sum^{m}_{i=1} w(x_i) | f_n(x_i) - f(x_i) |^2 \bbI_{\Omega}(x_i)} 
+ 2 \sqrt{\frac1m \sum^{m}_{i=1} w(x_i) | e_i |^2 \bbI_{\Omega}(x_i) }
\\
& \leq 4 \sqrt{\frac1m \sum^{m}_{i=1} w(x_i) | f(x_i) - f_n(x_i) |^2} 
 + 4 \sqrt{\frac1m \sum^{m}_{i=1} w(x_i) | e_i |^2  }.
}
Here, in the first step we used the facts that $f_n - \hat{f} \in \cP$ and $E_1$ occurs. In the second step we used the triangle inequality and the fact that $y_i = f(x_i) + e_i$, and in the third step we used the fact that $\hat{f}$ is a solution of the least-squares problem \ef{wls-prob}. 
Now let $E_2$ be the event that
\bes{
\sqrt{\frac1m \sum^{m}_{i=1} w(x_i) | f(x_i) - f_n(x_i) |^2} \leq \sqrt{\frac{3}{\epsilon}} \nm{f - f_n}_{L^2_{\varrho}(D)}
}
and $E_3$ be the event that
\bes{
\sqrt{\frac1m \sum^{m}_{i=1} w(x_i) | e_i |^2  } \leq \sqrt{\frac{3}{\epsilon}} \nm{e}_2.
}
Hence
\be{
\label{L2-bound-sub-E123-partial}
\text{$E_1$, $E_2$, $E_3$ occur }\Rightarrow\ \nm{f_n - \hat{f}}_{L^2_{\varrho}(\Omega)} \leq \frac{8}{\sqrt{\epsilon}} \nm{f-f_n}_{L^2_{\varrho}(D)} + \frac{8}{\sqrt{\epsilon}} \nm{e}_2.
}
Combining this with \ef{32C-today}, we conclude the following:
\be{
\label{L2-bound-sub-E123}
\text{$E_1$, $E_2$, $E_3$ occur }\Rightarrow\ \nm{f - \hat{f}}_{L^2_{\varrho}(\Omega)}  \leq \frac{9}{\sqrt{\epsilon}} \nm{f - f_n}_{L^2_{\varrho}(D)} + \frac{8}{\sqrt{\epsilon}} \nm{e}_{2}.
}
 Similarly, for the uniform norm, we first write 
\eas{
\nm{f - \hat{f}}_{L^{\infty}_{\varrho}(\Omega)}  & \leq \nm{f - f_n}_{L^{\infty}_{\varrho}(\Omega)} + \nm{f_n - \hat{f}}_{L^{\infty}_{\varrho}(\Omega)} 
\\
& \leq \nm{f - f_n}_{L^{\infty}_{\varrho}(\Omega)} + \sqrt{\kappa}  \nm{f_n - \hat{f}}_{L^{2}_{\varrho}(\Omega)} . 
}
Using the previous inequality, we deduce that
\be{
\label{Linf-bound-sub-E123}
\begin{split}
&\text{$E_1$, $E_2$, $E_3$ occur }\Rightarrow
\\
&~~~~ \ \nm{f - \hat{f}}_{L^{\infty}_{\varrho}(\Omega)} \leq  \nm{f - f_n}_{L^{\infty}_{\varrho}(\Omega)} 
 + 8\sqrt{\frac{\kappa}{\epsilon}} \nm{f-f_n}_{L^2_{\varrho}(D)} + 8\sqrt{\frac{\kappa}{\epsilon}} \nm{e}_2.
\end{split}
}
It therefore remains to estimate the probability that the events $E_1,E_2,E_3$ occur.

The case of $E_1$ involves a standard argument using the matrix Chernoff bound. Consider $\cP$ as a subspace of $L^2_{\varrho}(\Omega)$ of dimension $1\leq \tilde{n} \leq n$. 
Note that the lower bound follows from the fact that $1 \in \cP$. Let $\{ p_i \}^{\tilde{n}}_{i=1}$ be an orthonormal basis and observe that we can write
\bes{
\alpha = \lambda_{\min}(G),\quad \text{where }
G = \left ( \frac1m \sum^{m}_{i=1} w(x_i) \bbI_{\Omega}(x_i) \overline{p_j(x_i)} p_k(x_i) \right )^{n'}_{j,k=1}. 
}
This follows directly Parseval's identity and the fact that $\lambda_{\min}(G) = \inf \{ c^* G c : c \in \bbC^{n'},\ \nm{c}_2 = 1 \}$ for any self-adjoint, nonnegative definite matrix $G$.
In turn, we can write $G = X_1 + \cdots + X_m$, where the $X_i$ are i.i.d. copies of the self-adjoint, nonnegative definite random matrix
\bes{
X =  \left ( \frac1m w(x) \bbI_{\Omega}(x) \overline{p_j(x)} p_k(x) \right )^{n'}_{j,k=1},
}
where $x \sim \mu$. Observe that $\bbE(G) = I$ by orthonormality of the $p_j$'s. It is also a short argument -- using the fact that $\lambda_{\max}(X) = \sup \{ c^* X c : c \in \bbC^{n'},\ \nm{c}_2 = 1 \}$ for any self-adjoint, nonegative definite matrix $X$ -- to show that 
\bes{
\lambda_{\max}(X) \leq \frac1m \esssup_{x \sim \varrho} \{ w(x) \cK(\cP , L^2_{\varrho}(\Omega) )(x) \} = : R.
}
The matrix Chernoff bound \cite[Thm.\ 1.1]{tropp2012user-friendly} now implies that
\bes{
\bbP ( \alpha \leq 1/2 ) \leq n \exp \left ( - \frac{\log(1/2)/2 + 1/2}{R} \right ) \leq n \exp \left ( - \frac{1}{8 R} \right ).
}
Due to the assumption on $m$, we deduce that $\bbP(E_1^c) \leq \epsilon / 3$.

We estimate the probabilities of $E_2$ and $E_3$ via Markov's inequality. For $E_2$, we observe that
\bes{
\bbE \frac1m \sum^{m}_{i=1} w(x_i) | f(x_i) - p(x_i) |^2 = \int_{D} w(x) |f(x)-p(x)|^2 \D \mu(x) = \nm{f - p}^2_{L^2_{\varrho}(D)}.
} 
Hence $\bbP(E^c_2) \leq \epsilon/3$.
Similarly, for $E_3$ we observe that
\bes{
\bbE \frac1m \sum^{m}_{i=1} w(x_i) | e_i|^2 = \sum^{m}_{i=1} |e_i|^2 \int_{D} w(x) \D \mu(x) = \nm{e}^2_2.
}
Therefore $\bbP(E^c_3) \leq \epsilon/3$ as well.

We have shown that $\bbP(E^c_i) \leq \epsilon /3$ for $i = 1,2,3$. By the union bound, $\bbP(E_1 \cap E_2 \cap E_3) \geq 1-\epsilon$. To complete the proof, we now recall \ef{L2-bound-sub-E123} and \ef{Linf-bound-sub-E123}. \qed
}

\subsection{Recovery in $D$}

We next consider recovery over the whole of $D$. For $1 \leq q,r \leq \infty$, define the constant
\be{
\label{Remez-const}
\cR(\cP , L^q_{\varrho}(\Omega) , L^r_{\varrho}(D) ) = \sup \left \{ \nm{p}_{L^q_{\varrho}(D)} : p \in \cP,\ 0 < \nm{p}_{L^r_{\varrho}(\Omega) } \leq 1 \right \}.
}
We refer to this as a \textit{Remez constant} for $\cP$. Remez constants are well-studied objects in approximation theory when $\cP$ is a space of algebraic or trigonometric polynomials. In the current setting, $\cP$ may be arbitrary, but we adopt the same name regardless.

\thm{
\label{thm:accuracy-fulldomain}
Consider the setup of the previous lemma. Then the following holds with probability at least $1-\epsilon$. For any $e \in \bbC^m$ and any minimizer $\hat{f}$ of \ef{wls-prob}, the error satisfies
\eas{
\nm{f - \hat{f}}_{L^2_{\varrho}(D)} &\leq \nm{f - f_n}_{L^2_{\varrho}(D)} + \frac{8 \cR(\cP , L^2_{\varrho}(\Omega) , L^2_{\varrho}(D) ) }{ \sqrt{\epsilon}} \left ( \nm{f - f_n}_{L^2_{\varrho}(D)} + \nm{e}_2 \right )
\\
\nm{f - \hat{f}}_{L^{\infty}_{\varrho}(D)} &\leq \nm{f - f_n}_{L^{\infty}_{\varrho}(D)} + \frac{8 \cR(\cP , L^{\infty}_{\varrho}(\Omega) , L^{2}_{\varrho}(D) ) }{ \sqrt{\epsilon} } \left ( \nm{f - f_n}_{L^2_{\varrho}(D)} + \nm{e}_2 \right ),
}
where $f_n$ is as in \ef{orth-proj}.
}
\prf{
By the triangle inequality, 
\bes{
\nm{f - \hat{f}}_{L^2_{\varrho}(D)} \leq \nm{f - f_n}_{L^2_{\varrho}(D)} + \nm{f_n - \hat{f}}_{L^2_{\varrho}(D)}
}
for any $p \in \cP$. We now apply the Remez constant to obtain
\bes{
\nm{f - \hat{f}}_{L^2_{\varrho}(D)} \leq \nm{f - f_n}_{L^2_{\varrho}(D)} + \cR(\cP , L^2_{\varrho}(\Omega) , L^2_{\varrho}(D) ) \nm{f_n - \hat{f}}_{L^2_{\varrho}(\Omega)}.
}
To conclude, we apply \ef{L2-bound-sub-E123-partial}. The case for the $L^{\infty}$-norm error is identical.
\qed
}

\subsection{General recipe}

Theorem \ref{thm:accuracy-fulldomain} provides a general recipe for selecting a good sampling measure $\mu$ in the presence of a nonuniform cost function $c$. Specifically, given an \textit{accuracy parameter} $C > 0$, in view of \ef{C-exp-def} and \ef{m-cond-subdomain}, we strive to choose $\mu \ll \varrho$ to make the \textit{cost bound}
\bes{
\int_{\mathrm{supp}(\mu)} \frac{c(x)}{w(x)} \D \varrho(x) \cdot \esssup_{x \sim \varrho|_{\Omega}} \left \{ w(x) \cK(\cP , L^2_{\varrho}(\Omega))(x) \right \},
}
as small as possible, where $w$ is as in \ef{weight-fn} and $\Omega$ is any measurable set in which $v = \D \mu / \D \varrho $ is positive almost everywhere and for which
\bes{
\cR(\cP , L^2_{\varrho}(\Omega) , L^2_{\varrho}(D) ) \leq C .
}
Notice that $\Omega$ is purely a theoretical consideration, while $\mu$ forms part of the resulting algorithm. In particular, it needs to be computationally feasible to sample from $\mu$. Also, $\mu$ may also be required to have certain problem-dependent features -- for example, one may wish to consider choices of $\mu$ that lend themselves to hierarchical approximation schemes. We consider this issue further in the next section, where we introduce two approaches for designing $\mu$ in the case of the main example.

\section{Recovery using polynomial estimators}\label{s:recov-poly}

We now focus on the main example (see \S \ref{ss:main-examp}). Later in this section, we propose two strategies for choosing $\mu$, in \S \ref{ss:poly-approach1} and \S \ref{ss:poly-approach2}, respectively.

\subsection{Estimates for the Remez constant}

Before doing so, we first require estimates for the Remez constant \ef{Remez-const}.

\lem{
\label{lem:Remez}
Consider the main example and let $\Omega = (-(1-\sigma),1-\sigma)$ for some $0 < \sigma < 1$. Then the Remez constant \ef{Remez-const} satisfies
\bes{
\cR(\cP , L^2_{\varrho}(\Omega) , L^2_{\varrho}(D) ) \leq \cR(\cP , L^{\infty}_{\varrho}(\Omega) , L^2_{\varrho}(D) ) \leq n (\rho_{\sigma})^{n-1},
}
where
\be{
\label{rho-sigma-def}
\rho_{\sigma} = \frac{1 + \sqrt{2 \sigma - \sigma^2}}{1-\sigma}.
}
}
\prf{
Since $\{ \phi_i \}^{n}_{i=1}$ is an orthonormal basis for $L^2_{\varrho}(D)$ the functions
\bes{
\Phi_i(x) = \frac{1}{\sqrt{1-\sigma}} \phi_i \left(\frac{x}{1-\sigma} \right ),\quad x \in \Omega,\ i = 1,\ldots,n
}
are an orthonormal basis of $L^2_{\varrho}(\Omega)$. Let $p \in \cP$ be arbitrary and write $p = \sum^{n}_{i=1} c_i \Phi_i$, where $\nm{c}_2 = \nm{p}_{L^2_{\varrho}(\Omega)}$. By the Cauchy--Schwarz inequality, we have
\bes{
\nm{p}_{L^2_{\varrho}(D)} \leq \nm{p}_{L^{\infty}_{\varrho}(D)} \leq \sum^{n}_{i=1} \nm{\Phi_{i}}_{L^{\infty}_{\varrho}(D)} | c_i | \leq \sqrt{\sum^{n}_{i=1} \nm{\Phi_i}^2_{L^{\infty}_{\varrho}(D)} } \nm{p}_{L^2_{\varrho}(\Omega)}.
}
A classical result of Bernstein (see, e.g., \cite[Lem.\ 3.2]{platte2011impossibility}) now implies that
\bes{
\nm{\Phi_i}_{L^{\infty}_{\varrho}(D)} = \nm{\phi_i}_{L^{\infty}(-1/(1-\sigma),1/(1-\sigma))} \leq (\rho_{\sigma})^{i-1} \nm{\phi_i}_{L^{\infty}(-1,1)} = \sqrt{2i-1} (\rho_{\sigma})^{i-1},
}
where $\rho_{\sigma}$ is the unique solution to $(\rho_{\sigma} + 1/\rho_{\sigma})/2 = 1/(1-\sigma)$ in $(1,\infty)$, which is exactly by \ef{rho-sigma-def}. This gives
\bes{
\nm{p}_{L^2_{\varrho}(D)} \leq \nm{p}_{L^{\infty}_{\varrho}(D)} \leq \sqrt{\sum^{n}_{i=1} (2i-1) (\rho_{\sigma})^{2(i-1)} } \nm{p}_{L^2_{\varrho}(\Omega)} \leq n (\rho_{\sigma})^{n-1}  \nm{p}_{L^2_{\varrho}(\Omega)},
}
which completes proof. \qed
}

Combining this lemma with Theorem \ref{thm:accuracy-fulldomain} has the following consequence.

\cor{
\label{cor:main-examp-recov}
Consider the main example, where $D = (-1,1)$, $\varrho$ is the uniform probability measure and $\cP = \bbP_{n-1}$. Let $r > 1$, $n \geq  \log(r)/\log(1+\sqrt{8})$ and 
\be{
\label{sigma-scale-n}
0 \leq \sigma \leq (r^{1/n}-1)^2 / 16,
}
and suppose that $v = \D \mu / \D \varrho$ is positive almost everywhere in $(-(1-\sigma) , 1-\sigma)$.
Let $f \in L^2_{\varrho}(D)$ and $0 < \epsilon < 1$.
Then it is possible to construct a weighted least-squares estimator using $m \geq n$ samples $x_1,\ldots,x_m \sim_{\mathrm{i.i.d.}} \mu$ with an expected evaluation budget
\be{
\label{C-exp-bd-poly}
C_{\mathsf{exp}} \leq 8 \int_{\supp(\mu)} \frac{c(x)}{w(x)} \D x \cdot \sup_{|x| < 1-\sigma} \left \{ w(x) \cK(\cP , L^2_{\varrho}(-(1-\sigma),1-\sigma )(x) \right \} \cdot \log(3n/\epsilon),
}
where $w$ is as in \ef{weight-fn},
such that the following holds with probability at least $1-\epsilon$. For any $e \in \bbC^m$ and any minimizer $\hat{f}$ of \ef{wls-prob}, the error satisfies
\eas{
\nm{f - \hat{f}}_{L^2_{\varrho}(D)} & \leq \nm{f - f_n}_{L^2_{\varrho}(D)} + \frac{8 r n }{ \sqrt{\epsilon}} \left ( \nm{f - f_n}_{L^2_{\varrho}(D)} + \nm{e}_2 \right )
\\
\nm{f - \hat{f}}_{L^{\infty}_{\varrho}(D)} &\leq \nm{f - f_n}_{L^{\infty}_{\varrho}(D)} + \frac{8 r n }{ \sqrt{\epsilon}} \left ( \nm{f - f_n}_{L^2_{\varrho}(D)} + \nm{e}_2 \right ),
}
where $f_n$ is as in \ef{orth-proj}.
} 
\prf{
Let $\rho_{\sigma}$ be as in \ef{rho-sigma-def}. Then $\rho_{\sigma} \leq 1 + 4 \sqrt{\sigma}$ whenever $0 < \sigma \leq 1/2$. This holds in this case, due to the choice of $\sigma$ \ef{sigma-scale-n} and the assumption on $n$. Hence
\bes{
(\rho_{\sigma} )^n \leq (1 + 4 \sqrt{\sigma} )^n = (r^{1/n})^n = r.
}
The result now follows immediately from Theorem \ref{thm:accuracy-fulldomain}, Lemma \ref{lem:Remez}  and \ef{C-exp-def}. \qed
}

\subsection{First approach}\label{ss:poly-approach1}

In view of this result, we seek a sampling measure $\mu$ that minimizes the right-hand side of \ef{C-exp-bd-poly}. In our first approach, we assume the asymptotic growth of the cost function $c$ is known.
In light of the discussion in \S \ref{ss:main-examp}, we focus on cost functions that blow up at the endpoints $x = \pm 1$, as this is expected to be the most challenging case for polynomial approximation. Specifically, for some known $\alpha \geq 0$, we now assume that
\be{
\label{c-alg-growth}
c(x) \lesssim (1-x^2)^{-\alpha},\quad x \in (-1,1).
}

\thm{
\label{thm:first-approach}
Consider the setup of Corollary \ref{cor:main-examp-recov}, with $c$ as in \ef{c-alg-growth} for some $\alpha \geq 0$ and $r = 2$. Define the sampling measure
\be{
\label{mu-choice-1}
\D \mu(x) =  \frac{(1-x^2)^{\beta} \D x}{\int^{1}_{-1} (1-x^2)^{\beta} \D x},
}
where $\beta = - 1/2$ for $0 \leq \alpha < 1/2$ and $\beta = \alpha - 1+\delta$ for $\alpha \geq 1/2$, where $\delta > 0$ is arbitrary. Then the evaluation budget satisfies
\bes{
C_{\mathsf{exp}} \begin{cases} \lesssim_{\alpha} n \cdot \log(3n/\epsilon) & 0 \leq \alpha < 1/2 \\ \lesssim_{\alpha,\delta} n^{2(\alpha+\delta)} \cdot \log(3n/\epsilon) & \alpha \geq 1/2 \end{cases}
}
and the following holds with probability at least $1-\epsilon$. For any $e \in \bbC^m$, the weighted least-squares estimator $\hat{f}$ is unique and satisfies
\eas{
\nm{f - \hat{f}}_{L^2_{\varrho}(D)} &\lesssim \frac{n }{ \sqrt{\epsilon}} \left ( \nm{f - f_n}_{L^2_{\varrho}(D)} + \nm{e}_2 \right )
\\
\nm{f - \hat{f}}_{L^{\infty}_{\varrho}(D)} &\lesssim \nm{f - f_n}_{L^{\infty}_{\varrho}(D)} +  \frac{n }{ \sqrt{\epsilon}} \left ( \nm{f - f_n}_{L^2_{\varrho}(D)} + \nm{e}_2 \right ),
}
where $f_n$ is as in \ef{orth-proj}.
}
\prf{
Uniqueness of $\hat{f}$ follows from the fact that the $m \geq n$ sample points are distinct with probability one. Consequently, no non-zero polynomial $p \in \cP$ can vanish at the sample points, with probability one.
For the rest of the proof, we use Corollary \ref{cor:main-examp-recov}. Suppose first that $-1/2 \leq \beta \leq 0$ (the lower bound comes from the definition of $\beta$). Then we set $\sigma = 0$ in \ef{C-exp-bd-poly}. Notice that this choice of $\sigma$ trivially satisfies \ef{sigma-scale-n}. Since $w(x)^{-1} = (1-x^2)^{\beta} / \int^{1}_{-1} (1-x^2)^{\beta} \D x$ we get
\bes{
C_{\mathsf{exp}} \lesssim_{\alpha,\beta}  \sup_{|x| < 1} \{ (1-x^2)^{-\beta} \cK(\cP , L^2_{\varrho}(-1,1) )(x) \} \cdot \log(3n/\epsilon) .
}
Now \ef{Christoffel-UB} implies that $\cK(\cP , L^2_{\varrho}(-1,1))(x) \leq n^{2(1+\beta)} (1-x^2)^{\beta}$. Therefore
\bes{
C_{\mathsf{exp}} \lesssim_{\alpha,\beta} n^{2(1+\beta)} \cdot \log(3n / \epsilon).
}
This gives the result for $-1/2 \leq \beta \leq 0$.

Now let $\beta > 0$.
Consider the Christoffel function in \ef{C-exp-bd-poly}. Let $p \in \cP = \bbP_{n-1}$ and, for $x \in (-(1-\sigma),1-\sigma)$, write $p(x) = P(y)$, where $x = (1-\sigma) y$ for $y \in (-1,1)$. 
Then
\bes{
|p(x) |^2 = |P(y)|^2 \leq n^2 \nm{P}^2_{L^2_{\varrho}(D)} = \frac{n^2}{1-\sigma} \nm{p}^2_{L^2_{\varrho}(-(1-\sigma),1-\sigma))}.
}
Here, the factor $n^2$ follows from \ef{Christoffel-UB}.
Therefore
\bes{
\cK(\cP , L^2_{\varrho}(-(1-\sigma),1-\sigma))(x) \leq \frac{n^2}{1-\sigma}.
}
Using the definition of $\mu$, we deduce that
\bes{
C_{\mathsf{exp}} \lesssim_{\alpha,\beta} \frac{n^2}{1-\sigma} \sup_{|x| \leq 1-\sigma} \{  (1-x^2)^{-\beta}   \} \cdot \log(3n/\epsilon).
}
Now observe that $1 \geq 1-x^2 \geq 2 \sigma - \sigma^2 \geq \sigma$ as $0 < \sigma < 1$, and therefore
\bes{
\sup_{|x| \leq 1-\sigma} \{ (1-x^2)^{-\beta} \} \leq \sigma^{-\beta}.
}
Hence, if $0 < \sigma \leq 1/2$, we deduce that $C_{\mathsf{exp}} \lesssim_{\delta} n^2 \cdot \sigma^{-\beta} \cdot \log(3n/\epsilon)$.
Finally, we choose $\sigma$ maximally such that \ef{sigma-scale-n} holds with $r = 2$, i.e., $\sigma = \frac{(2^{1/n}-1)^2}{16}$. We then notice that
\be{
\label{sigma-asymp}
\sigma \sim \log(2)^2 / (16 n^2) + \ord{1/n^3},\quad n \rightarrow \infty.
}
This yields $C_{\mathsf{exp}} \lesssim_{\alpha,\beta} n^{2(1+\beta)} \cdot \log(3n/\epsilon)$, which gives the result for $\beta > 0$ as well.
\qed
}

This theorem shows that the choice \ef{mu-choice-1} leads to an expected evaluation cost that is log-algebraic in $n$. If $0 \leq \alpha < 1/2$, then we obtain near-optimal, log-linear cost in $n$. This is unsurprising, since the cost function \ef{c-alg-growth} is integrable with respect to the Chebyshev measure \ef{arcsine-meas} in this case. The conclusion of log-linear cost is, in this case, a straightforward consequence of \ef{Christoffel-UB} and Theorem \ref{t:err-bd-prob}. Conversely, when $\alpha \geq 1/2$ we obtain a cost that scales like $n^{2(\alpha+\delta)} \cdot \log(n)$, where $\delta > 0$ can be arbitrarily small. As we discuss below, this scaling is optimal up to $\delta$.

\subsection{Second approach}\label{ss:poly-approach2}

The disadvantage of the previous approach is that requires knowledge of the cost function and, for $\alpha \geq 1/2$, it leads to a sample complexity bound that is not quite optimal, due to the factor $\delta > 0$. We now consider a different approach that overcomes these issues.

\thm{
\label{thm:second-approach}
Consider the setup of Corollary \ref{cor:main-examp-recov}, with $c$ as in \ef{c-alg-growth} for some $\alpha \geq 0$ and $r = 2$. Define the sampling measure
\be{
\label{mu-sigma-adjust}
\D \mu(x) = ((1-\sigma) \pi)^{-1} \left ( 1 - (x/(1-\sigma))^2 \right )^{-1/2} \bbI_{(-1+\sigma,1-\sigma)}(x) \D x,
}
where $\sigma = \sigma(n) = \frac{(2^{1/n}-1)^2}{16}$. Then the evaluation budget satisfies
\bes{
C_{\mathsf{exp}} \lesssim_{\alpha} n^{\max \{ 1,2\alpha  \}} \cdot \log(3n/\epsilon)
}
and the following holds with probability at least $1-\epsilon$. For any $e \in \bbC^m$, the weighted least-squares estimator $\hat{f}$ is unique and satisfies
\eas{
\nm{f - \hat{f}}_{L^2_{\varrho}(D)} &\lesssim \frac{n }{ \sqrt{\epsilon}} \left ( \nm{f - f_n}_{L^2_{\varrho}(D)} + \nm{e}_2 \right )
\\
\nm{f - \hat{f}}_{L^{\infty}_{\varrho}(D)} & \lesssim \nm{f - f_n}_{L^{\infty}_{\varrho}(D)} +  \frac{n }{ \sqrt{\epsilon}} \left ( \nm{f - f_n}_{L^2_{\varrho}(D)} + \nm{e}_2 \right ) ,
}
where $f_n$ is as in \ef{orth-proj}.
}
\prf{
Uniqueness follows by the same argument as before. Now consider the Christoffel function appearing in \ef{C-exp-bd-poly}. As in the previous proof, let $p \in \cP = \bbP_{n-1}$ and write $p(x) = P(y)$, where $x = (1-\sigma) y$. Then $\nm{p}^2_{L^2_{\varrho}(\Omega)} = (1-\sigma) \nm{P}^2_{L^2_{\varrho}(D)}$ and therefore
\bes{
\cK(\cP , L^2_{\varrho}(\Omega))(x) = \frac{1}{1-\sigma} \cK(\cP , L^2_{\varrho}(D))(y) .
}
Hence \ef{Christoffel-UB} gives
\bes{
\cK(\cP , L^2_{\varrho}(\Omega))(x) \leq \frac{4 n}{\pi(1-\sigma)(1-(x/(1-\sigma))^2)^{1/2} },\quad |x| < 1-\sigma.
}
Using the definition of $w$, we obtain
\bes{
\sup_{|x| \leq 1-\sigma} \{ w(x) \cK(\cP , L^2_{\varrho}(\Omega))(x) \} \lesssim n,
}
and plugging this into  \ef{C-exp-bd-poly} we deduce that
\bes{
C_{\mathsf{exp}} \lesssim \int^{1}_{-1} \frac{c(x)}{w(x)} \D x \cdot n \cdot \log(3n/\epsilon). 
}
Consider the integral. Using the definitions of $c$ and $w$, we have
\eas{
\int^{1}_{-1} \frac{c(x)}{w(x)} \D x & \lesssim \int^{1-\sigma}_{-1+\sigma} (1-x^2)^{-\alpha} \left ( 1 - (x/(1-\sigma))^2 \right )^{-1/2} \D x
\\
& \lesssim \int^{1-\sigma}_{0} (1-x)^{-\alpha} \left (1 - x/(1-\sigma) \right )^{-1/2} \D x = :  I .
}
We now make a change of variables and integrate by parts, to get
\eas{
I & = (1-\sigma)^{1/2} \int^{1}_{\sigma} y^{-\alpha} (y - \sigma)^{-1/2} \D y
\\
& = 2 (1-\sigma) + 2(1-\sigma)^{1/2} \alpha \int^{1}_{\sigma} y^{-\alpha-1} (y-\sigma)^{1/2} \D y 
\\
& \leq 2 + 2\alpha \int^{1}_{\sigma} y^{-\alpha-1/2} \D y
\\
& = 2 + \frac{2\alpha}{1/2-\alpha} \left ( 1 - \sigma^{1/2-\alpha} \right ).
}
Hence $I \lesssim_{\alpha} \max \{ 1, \sigma^{1/2-\alpha}  \}$. We deduce that $C_{\mathsf{exp}} \lesssim_{\alpha} \max \{ 1, \sigma^{1/2-\alpha}  \} \cdot n \cdot \log(3n/\epsilon)$.
The result now follows from \ef{sigma-asymp}. \qed
}

The sampling measure considered in this theorem overcomes the limitation of that of Theorem \ref{thm:first-approach}, since it is independent of the cost function. The strategy involves sampling according to the Chebyshev measure on a subdomain -- this being a near-optimal sampling measure for polynomial approximation on a subdomain -- where the size of the domain is chosen judiciously to ensure a small Remez constant. The downside of this approach is that it is inherently \textit{nonhierarchical}. The sampling measure depends on $\sigma$, and therefore $n$. Hence, if one were to increase the polynomial degree (e.g., to compute a more accurate approximation) any existing samples could not be re-used, as they would be drawn from the wrong measure. Conversely, the approach developed in \S \ref{ss:poly-approach1} is trivially hierarchical, since the sampling measure is independent of $n$.

We remark in passing that standard Christoffel sampling is also nonhierarchical (see \S \ref{ss:CS}), as the measure \ef{mu-opt} also depends on $n$. However, there are hierarchical variants on Christoffel sampling \cite{adcock2020near-optimal,arras2019sequential,migliorati2019adaptive}. Whether these could be adapted to the current setting of nonuniform evaluation costs is an open problem.

\rem{
In both theorems, $\sigma$ is chosen so that the term $(\rho_{\sigma})^n \leq 2$ for all $n$. This leads to the behaviour \ef{sigma-asymp}, where $\sigma$ decays like $1/n^2$ as $n \rightarrow \infty$. It is precisely this decay that contributes to the higher evaluation cost scaling with respect to $n$. For instance, in the proof of Theorem \ref{thm:second-approach} we see that
\bes{
C_{\mathsf{exp}} \leq m \cdot \max \{ 1, \sigma^{1/2-\alpha} \}.
}
This raises the question of whether $\sigma$ can be chosen in a different way. It is readily verified that insisting on algebraic growth, i.e., $(\rho_{\sigma})^n \leq n^{\beta}$ for some $\beta > 0$, can at best improve the scaling by a polylogarithmic factor in $n$. One could instead demand that $\sigma$ decays at most algebraically in $n$. However, if $\sigma^{-1} = c n^{\beta}$ for some $0 < \beta < 2$ then $(\rho_{\sigma} )^n \gtrsim \exp(c' n^{1-\beta/2})$ for $c'$ depending on $c$. Thus, any attempt to reduce the algebraic scaling of $C_{\mathsf{exp}}$ below $2 \alpha$ necessarily results in an exponentially-growing term in the error bound.
}

\rem{
The algebraic scaling $n^{2 \alpha}$ for $\alpha \geq 1/2$ is effectively optimal. This follows from a recent result of \cite{xu2024stability}, which is based on earlier results in \cite{adcock2019optimal,platte2011impossibility}. We now briefly argue why this applies. Note first that the expected cost is only finite if and only if the density $v$ of the measure $\mu$ is such that $\int^{1}_{-1} (1-x^2)^{-\alpha} v(x) \D x < \infty$. In this case, $C_{\mathsf{exp}}$ is proportional to $m$. Assume, for convenience, that $v$ is symmetric an increasing on $[0,1]$. Then $v(x) \lesssim_{\alpha} (1-x^2)^{\alpha-1}$. This means that $\mu$ is a \textit{modified Jacobi} measure, with $c_1 = 0$ (see \cite{xu2024stability}). Therefore, \cite[Thm.\ 1.4]{xu2024stability} shows that the least-squares estimator is exponentially ill-conditioned if $m$ scales like $n^{\beta}$ for \textit{any} $\beta < 2\alpha$. Theorems \ref{thm:first-approach} and \ref{thm:second-approach} both imply the estimator is well-conditioned. Hence, they achieve the optimal scaling.

Moreover, a stronger result also holds. Since the estimators in Theorems \ref{thm:first-approach} and \ref{thm:second-approach} yield quasi-best polynomial approximations, they must converge exponentially fast in $n \approx m^{1/(2\alpha)}$ whenever $f$ is analytic. The result \cite[Thm.\ 1.5]{xu2024stability} shows that no well-conditioned estimator can achieve a faster rate of exponential convergence in $m$.
}

\subsection{Numerical example}
\begin{figure}[t]
\begin{center}
\begin{tabular}{cc}
\includegraphics[width = 0.45\textwidth]{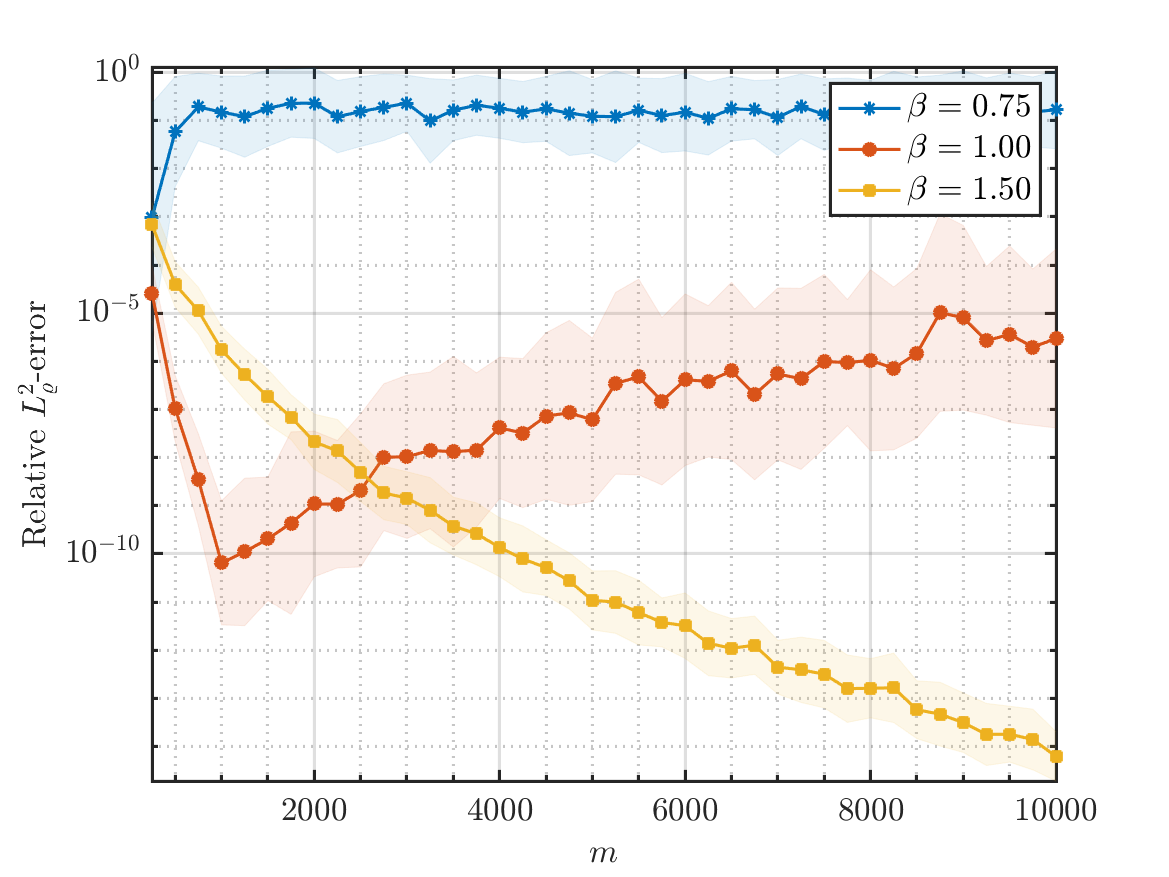} &
\includegraphics[width = 0.45\textwidth]{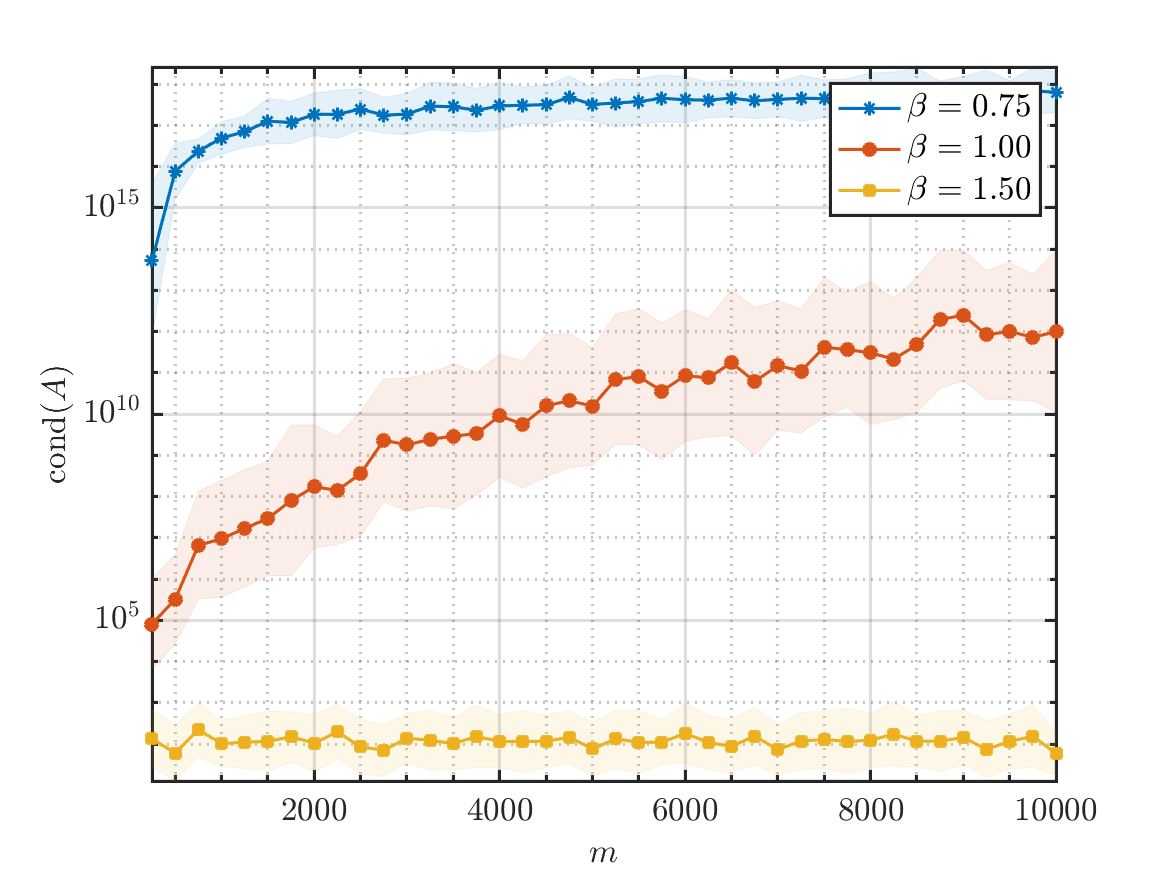}
\end{tabular}
\end{center}
\caption{Error (left) and condition number (right) for the least-squares estimator based on the sampling measure $\D \mu(x) \propto (1-x^2)^{\alpha-1} \D x$ for $\alpha = 1.5$. We use the scaling $m = 0.5 \cdot n^{2 \beta}$ for various $\beta$. The function $f(x) = (1.1-x)^{-1}$. We use 50 trials for each value of $m$. The solid line is the (geometric) mean and the shaded region shows one (geometric) standard deviation. }
\label{fig1}
\end{figure}

We conclude with a numerical example. In Fig.\ \ref{fig1} we plot the error (left) and condition number (right) for the least-squares polynomial estimator with cost function with $\alpha = 3/2$. When $m = n^{1.5}$, the estimator is extremely ill-conditioned and does not converge as $m \rightarrow \infty$ (see the previous remark). When $m = n^2$ the estimator initially reaches a small error, but the error drifts as $m$ increases further, indicating an ill-conditioned estimator. Conversely, when $m = n^3$, we obtain a well-conditioned, convergent estimator, as predicted by our theoretical results.

\begin{SCfigure}[][t!]
\includegraphics[width = 0.45\textwidth]{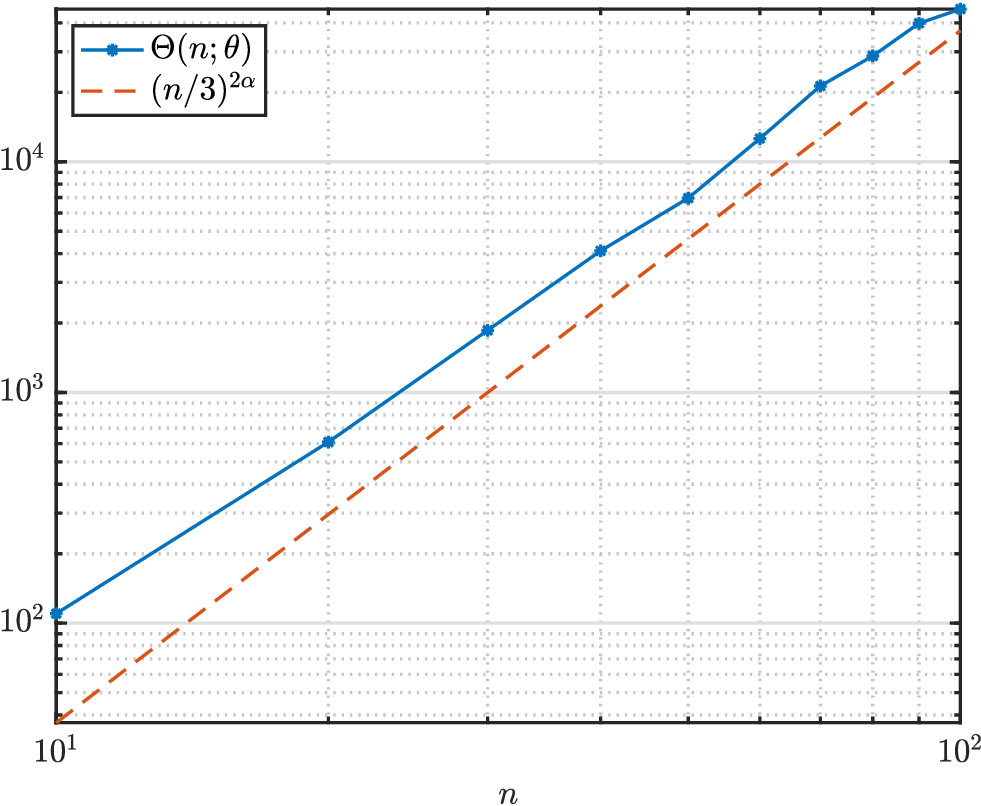}
\caption{The quantity $\Theta(n ; \theta)$ versus $n$, where $\Theta(n ; \theta)$ is the smallest value of $m$ such that the mean condition number of the least-squares estimator is at most $\theta$. In this experiment, we consider $\alpha = 1.5$ and $\theta = 10$. We estimate $\Theta(n ; \theta)$ using 50 trials, incrementing $m$ from $m = n$ in steps of 50 until the estimated mean is at most $\theta$.}
\label{fig2}
\end{SCfigure}

Next, in Figure \ref{fig2} we empirically verify the scaling $n^{2 \alpha}$. We do this by finding, for each $n$, the minimal value of $m$ such that the least-squares estimator is well-conditioned. This figure shows good agreement between the empirically-estimated scaling and the theoretically-optimal scaling $n^{2 \alpha }$.

\section{Conclusions and open problems}\label{s:conclusions}

This paper presents a first foray into function recovery with nonuniform evaluation costs -- a topic that, to the best of the author's knowledge, has not been widely studied. Focusing on linear least-squares estimators, we derived a general recipe for accurate recovery, which involves the Christoffel function and Remez constant of the subspace $\cP$. We then applied this to the case of polynomial approximation, yielding two methods which were shown to be near-optimal for algebraically-growing cost functions.

There are many avenues for further investigation. Two specific open problems are there the extension of the main example to the multivariate case and the question of designing hierarchical sampling measures. More broadly, there are questions about the information complexity of recovery in different function classes. These questions are well-studied in the case of uniform costs \cite{krieg2021function,krieg2021functionII}. However, to the best of the author's knowledge, little is known about the case of nonuniform evaluation costs.

\section*{Acknowledgements}

BA acknowledges the support of NSERC through grant RGPIN-2021-611675.  The idea for this paper arose during discussions at MCQMC 2024. He would like to thank the organizers and participants for a providing stimulating research environment.

\bibliographystyle{spmpsci}
\bibliography{refs}

\end{document}